\documentclass[12pt]{article}       

\usepackage{amsmath,amsthm}
\usepackage{amssymb}
\usepackage{latexsym}
\begin{document}

\newcommand{\comment}[1]{}    
\newcommand{\hs}{\enspace}
\newcommand{\hhs}{\thinspace}
\newcommand{\real}{\ifmmode {\rm R} \else ${\rm R}$ \fi}

\newtheorem{theorem}{Theorem}
\newtheorem{lemma}[theorem]{Lemma}         
\newtheorem{corollary}[theorem]{Corollary}
\newtheorem{definition}[theorem]{Definition}
\newtheorem{claim}[theorem]{Claim}
\newtheorem{conjecture}[theorem]{Conjecture}
\newtheorem{proposition}[theorem]{Proposition}

\newtheorem{remark}[theorem]{Remark}



\def\edge{\leftrightarrow}
\def\noedge{\not\leftrightarrow}
\def\twoedge{\Leftrightarrow}
\def\to{\rightarrow}
\def\Hrl{H^{(r)}_{l+1}}
\def\Krl{{\cal K}^{(r)}_l}
\def\Krl{{\cal K}^{(r)}_{l+1}}
\def\cF{{\cal F}}
\def\cG{{\cal G}}
\def\cH{{\cal H}}
\def\e{\varepsilon}
\def\bF{
{\cal {\bf F}}}
\def\odel{o_{\delta}}
\def\od1{o_{\delta}(1)}


\title{\bf Counting substructures III: quadruple systems}
\author{Dhruv Mubayi
\thanks{ Department of Mathematics, Statistics, and Computer
Science, University of Illinois, Chicago, IL 60607.  email:
mubayi@math.uic.edu; research  supported in part by  NSF grant DMS
0653946
\hfil\break\null\hskip .23in Keywords: {\it Hypergraph Tur\'an
numbers, removal lemma, stability theorems}}}
\date{\today}
\maketitle

\begin{abstract}
For various quadruple systems $F$, we give asymptotically sharp lower bounds on the number
of copies of $F$ in  a quadruple system with a prescribed number of
vertices and edges.  Our results extend those of F\"uredi, Keevash, Pikhurko, Simonovits and  Sudakov who proved under the same conditions that there is one copy of $F$.
Our proofs use the hypergraph removal Lemma and
stability results for the corresponding Tur\'an problem proved by the above authors.

\end{abstract}

\section{Introduction}

Given a $k$-uniform hypergraph ($k$-graph for short) $F$, let ex$(n, F)$, the Tur\'an number of $F$, be the maximum number of edges in an $n$ vertex $k$-graph with no copy of $F$.  Beginning with a result of Rademacher the following phenomenon has been discovered and studied intensely for $k=2$:  If an $n$  vertex $k$-graph has
ex$(n, F)+1$ edges, then it has not just one, but many copies of $F$. This is the third in a series of papers where we refine and extend earlier results of Rademacher, Erd\H os \cite{E1, E2}, Lovasz-Simonovits \cite{LS} and others on this phenomenon. The first
paper in this series \cite{M1} studied the question when $F$ is a color critical graph, the second \cite{M2} studied the case of 3-graphs, and here we study 4-graphs.

 We are able to go beyond earlier work in this area  due to the new tool we have at our disposal: the hypergraph removal lemma.  This is a consequence of the hypergraph regularity lemma proved by Gowers \cite{G}, Nagle-R\"odl-Schacht \cite{NRS}, R\"odl-Skokan \cite{RS}, Tao \cite{T}.

\begin{theorem} {\bf (Hypergraph Removal Lemma \cite{G, NRS, RS, T})} \label{removal}
Fix $k \ge 2$ and a $k$-graph $F$ with $f$ vertices. Suppose that an
$n$ vertex $k$-graph $\cH$ has at most $o(n^f)$ copies of $F$. Then
there is a set of edges in $\cH$ of size $o(n^k)$ whose removal from
$\cH$ results in a $k$-graph with no copies of $F$.
\end{theorem}

There are two types of configurations we will be concerned with in this paper, books and expanded triangles.  These (with one exception that was treated in \cite{M2}) represent all known cases of a 4-graph $F$ where ex$(n,F)$ has been exactly determined.

 For $2 \le l \le 4$, the $l$-book  $P_l$ is the 4-graph with $l+1$ edges $l$ of which  share the same three points, and another edge that contains the remaining point in each of the $l$ edges together with $4-l$ new points. Explicitly
\begin{align} P_2&=\{123a, 123b, abcd\} \notag \\
P_3&=\{123a, 123b, 123c, abcd\}\notag \\
P_4&=\{123a, 123b, 123c, 123d, abcd\}\notag \end{align}
The expanded triangle $C_3$ is the 4-graph obtained from a graph triangle by replacing each vertex by a pair of vertices. Formally,
$$C_3=\{1234, 3456, 1256\}.$$

\begin{definition}
Suppose  $F$ is a 4-graph with the property that for sufficiently large $n$, there is a unique (up to isomorphism) 4-graph $\cH(n, F)$ with ex$(n,F)$ edges. Let $c(n,F)$ be the minimum number of copies of $F$ in the 4-graph obtained from $\cH(n, F)$ by adding an edge, where the minimum is taken over all possible edges that may be added.
 \end{definition}

  In the next two subsections will discuss $\cH(n,F)$ and $c(n, F)$ for $F \in \{P_2, P_3, P_4, C_3\}$ and then we will state our results.

Notation: We associate a hypergraph with its edge set. The number of edges in a hypergraph $\cH$ is $|\cH|$. Given hypergraphs
$F, \cH$ ($F$ has $f$ vertices), a copy of $F$ in $\cH$ is a subset of
$f$ vertices and $|F|$ edges of $\cH$ such that the subhypergraph formed
by this set of vertices and edges is isomorphic to $F$.    In other words, if we denote $Aut(F)$ to be the number of automorphisms of $F$, then the number of copies of $F$ in $\cH$ is the number of edge-preserving injections from $V(F)$ to $V(\cH)$ divided by $Aut(F)$.
For a set $S$ of vertices, define $d_{\cal H}(S)$ to be the number of edges of $\cH$ containing $S$.  If $S=\{v\}$, we simply write $d_{\cal H}(v)$.
We will omit floor and ceiling symbols whenever they are not crucial, so that the presentation is clearer.

\subsection{Books}
{\bf $P_2=\{123a, 123b, abcd\}$.}
Write
$$t^4(n)=\left\lfloor {n\over 4} \right\rfloor \left\lfloor {n+1\over 4} \right\rfloor \left\lfloor {n+2\over 4}
\right \rfloor\left\lfloor {n+3\over 4} \right\rfloor$$ for the number of edges in $T^4(n)$, the complete
4-partite 4-graph with the maximum number of edges.  It is easy to see that $T^4(n)$ contains no copy of $P_2$.
Frankl and F\"uredi  \cite{FF} conjectured, and Pikhurko \cite{P} proved, that
ex$(n, P_2)=t^4(n)$ for $n$ sufficiently large.  This shows that $c(n, P_2)$ is defined and one achieves $c(n, P_2)$ by adding an edge to $T^4(n)$ with two points in each of two parts, and no point in the remaining two parts.  We then see that
$$c(n, P_2)=2(n/4)^3-O(n^2) =\Theta(n^3).$$

\medskip
{\bf $P_3=\{123a, 123b, 123c, abcd\}$.}
Say that a 4-graph has a $(2,2)$-partition if it has a vertex partition into two parts so that every edge intersects each part in two points.
Write
$$d^4(n)={\lfloor n/2\rfloor \choose 2}{\lceil n/2 \rceil \choose 2}$$ for the number of edges in $D^4(n)$, the $n$ vertex 4-graph with a $(2,2)$-partition having
the maximum number of edges.  It is easy to see that $D^4(n)$ contains no copy of $P_3$.
F\"uredi, Simonovits and Pikhurko \cite{P} proved that
ex$(n, P_3)=d^4(n)$ for $n$ sufficiently large.   This shows that $c(n, P_3)$ is defined and one achieves $c(n, P_2)$ by adding an edge to $D^4(n)$ with exactly three points in the part of size $\lceil n/2 \rceil$.  We then see that
$$c(n, P_3)=4{\lfloor n/2\rfloor -1 \choose 2}(\lceil n/2 \rceil -3)
=2(n/2)^3-O(n^2) =\Theta(n^3).$$

\medskip
{\bf $P_4=\{123a, 123b, 123c, 123c, abcd\}$.}
 A 4-graph $\cH$ is odd if it has
a vertex partition $A \cup B$ such that every edge intersects both
parts in an odd number of vertices.  Let $B^4(n)$ be the odd 4-graph with the maximum
number of edges. Note that
$$b^4(n):=|B^4(n)|=\max_{1 \le a \le n}
 {a\choose 3}(n-a) +{n-a \choose 3}a$$
  is {\em not} achieved by choosing
$a=\lfloor n/2 \rfloor$, but it can easily be shown that $|a-n/2|<\sqrt{3n}/2+1$.
F\"uredi, Mubayi and Pikhurko \cite{P} proved that
ex$(n, P_4)=b^4(n)$ for $n$ sufficiently large.  This shows that $c(n, P_4)$ is defined and one achieves $c(n, P_2)$ by adding an edge to $D^4(n)$ with two points in each part. We then see that
$$c(n,P_4)=4{n/2 \choose 3}-O(n^2)=\Theta(n^3).$$

Now we state our result about counting books.

\begin{theorem} \label{books} Fix $l \in \{2,3,4\}$.
For every $\e>0$ there exists $\delta>0$ and $n_0$ such that the
following holds for $n>n_0$. Let $\cH$ be an $n$ vertex  4-graph with ex$(n, P_l)+q$
edges where $q<\delta n$.  Then the number of copies of $P_l$ in
$\cH$ is at least $q(1-\e)c(n, P_l)$. The expression $q$ is sharp for
$1\le q<\delta n$. Moreover, if the number of copies is less than
$\delta n^4$, then there is a collection of $q$ distinct edges that
each lie in $(1-\e)c(n, P_l)$ copies of $P_l$ with no two of these edges accounting for the same copy of $P_l$.
\end{theorem}

\subsection{Expanded triangle}

 The expanded triangle is a 4-uniform example whose extremal value has been studied by Frankl \cite{F}, Sidorenko, Keevash and Sudakov \cite{KSFr}.
Recall that $C_3=\{1234, 3456, 1256\}$ is the 4-graph obtained from a graph triangle by expanding each vertex to a set of size two.
Frankl \cite{F} proved that every $n$ vertex 4-graph containing no copy of  $C_3$ has at most $(1+o(1))b^4(n)$ edges. Recently, Keevash and Sudakov \cite{KSFr} sharpened this by proving that the unique 4-graph that achieves this maximum is $B^4(n)$.
Adding an edge to $B^4(n)$ results in at least 
$$c(n,C_3)=3(n/2)^2-O(n)=\Theta(n^2)$$
 copies of $C_3$.

\begin{theorem} \label {C3q}
For every $\e>0$ there exists $\delta>0$ and $n_0$ such that the
following holds for $n>n_0$. Let $\cH$ be an $n$ vertex  4-graph with $b^4(n)+q$
edges where $q<\delta n^2$.  Then the number of copies of $C_3$ in
$\cH$ is at least $q(1-\e)c(n, C_3)$. The expression $q$ is sharp for
$1\le q<\delta n^2$. Moreover, if the number of copies is less than
$\delta n^4$, then there is a collection of $q$ distinct edges that
each lie in $(1-\e)c(n, C_3)$ copies of $C_3$ with no two of these edges accounting for the same copy of $C_3$.
\end{theorem}

We remark that although our proof follows the same general structure as that in  \cite{KSFr}, some new ideas are needed.  In particular, since we start our proof with an application of the removal lemma, we do not have such fine control over the size of the parts in the underlying hypergraph as in \cite{KSFr}.  Thus our approach is somewhat more robust, although the approach in \cite{KSFr} extends to the $k$-uniform case which we do not address here.

Throughout the paper we will frequently use the notation $\delta \ll \e$, which is supposed to mean that $\delta$, and any function of $\delta$ (that tends to zero with $\delta$) used in a proof is smaller than any function of $\e$ used in the proof.  It is pretty difficult to write the precise dependence between $\delta$ and $\e$ as one of the constraints comes from an application of the removal lemma.

\section{Counting  $P_2$'s}

Recall that $c(n,P_2)=2(n/4)^3+\Theta(n^2)$. Theorem \ref{books} for $l=2$  follows from the following result.

\begin{theorem} \label{P2}
For every $\e>0$ there exists $\delta>0$ and $n_0$ such that the
following holds for $n>n_0$. Every $n$ vertex 4-graph with
$t^4(n)+1$ edges contains either

$\bullet$ an edge that lies in at least $(2-\e)(n/4)^3$ copies of
$P_2$, or

$\bullet$ at least $\delta n^4
$ copies of $P_2$.
\end{theorem}

\bigskip
\noindent {\bf Proof of Theorem \ref{books} for $l=2$.}
Remove $q-1$ edges
from $\cH$ and apply Theorem \ref{P2}.
If we find $\delta n^4$
copies of $P_2$, then since $q<\delta n$, the number of copies is
much larger than $(1-\e)qc(n, P_2)$ and we are done.
Consequently, we find
an edge $e_1$ in at least $(2-\e)(n/4)^3>(1-\e)c(n, P_2)$ copies of $P_3$. Now remove $q-2$ edges from $\cH-e_1$ and repeat this argument to obtain $e_2$.  In this way we obtain edges $e_1, \ldots, e_q$ as required.

The bound is sharp due to the following construction.  Suppose that
$T^4(n)$ has parts $W, X, Y, Z$.   Fix a pair of points $a,b \in W$, and add $q$ edges of the form $abcd$ where $c,d \in X$.
It is easy to see that each added edge lies in
$2(n/4)^3+O(n^2)$ copies of $P_2$ and no copy of $P_2$ contains two of
the new edges. Consequently, the copies of $P_2
$ are counted exactly
once. \qed

We will need the following stability theorem for $P_2$ proved by Pikhurko \cite{P}

\begin{theorem} {\bf ($P_2$ stability \cite{P})}\label{P2stability}
Let $\cH$ be a 4-graph with $n$ vertices and $t^4(n)-o(n^4)$ edges
that contains no copy of $P_2$. Then there is a partition of the
vertex set of $\cH$ into $W \cup X \cup Y \cup Z$ so that the number of edges that
intersect a part in at least two  points is $o(n^4)$. In other words, $\cH$ can
be obtained from $T^4(n)$ by adding and deleting a set of $o(n^4)$
edges.
\end{theorem}

\bigskip

\noindent{\bf Proof of Theorem \ref{P2}.}
 Given $\e$ let $0<\delta \ll \e$.
Write $o_{\delta}(1)$ for any function that approaches zero as
$\delta$ approaches zero and moreover, $o_{\delta}(1) \ll \e$. Let
$n$ be sufficiently large and let $\cH$ be an $n$ vertex 4-graph
with $t^4(n)+1$ edges.  Write $\#P_2$ for the number of copies of
$P_2$ in $\cH$.

We first argue that we may assume that $\cH$ has minimum degree at
least  $d=(3/32)(1-\delta_1){n\choose 3}$, where
$\delta_1=\delta^{1/4}$. Indeed, if this is not the case, then
remove a vertex of degree less than $d$ to form the 4-graph $\cH_1$
with $n-1$ vertices.  Continue removing a vertex of degree less than
$d$  if such a vertex exists. If we could continue this process for
$\delta_2 n$ steps, where $\delta_2=\delta^{1/2}$, then the
resulting 4-graph $\cH'$ has $(1-\delta_2)n$ vertices and number of
edges at least
\begin{align}
t^4(n)-\delta_2 n d &>{3\over 32}(1-\delta-4\delta_2(1-\delta_1)){n \choose 4}\notag \\
&={3\over 32}(1-\delta-4\delta_2+4\delta_1\delta_2){n \choose 4}\notag \\
&>{3\over 32}(1-\delta-4\delta_2+(2\delta+6\delta_2^2+6\delta_2^2\delta-4\delta_2\delta)){n \choose 4}\notag \\
&={3\over 32}(1+\delta)(1-4\delta_2+6\delta_2^2){n \choose 4} \notag \\
&>{3\over 32}(1+\delta)(1-\delta_2)^4{n \choose 4} \notag \\
&>{3\over
32}(1+\delta){(1-\delta_2)n \choose 4} \notag.
\end{align}
By the result of
Pikhurko \cite{P} and Erd\H os-Simonovits supersaturation
we conclude that $\cH$ has at least $\delta'n^7>\delta n^4$ copies of $P_2$
(for some fixed $\delta'>0$ depending on $\delta$) and we are done.  So we may assume that
this process of removing vertices of degree less than $d$ terminates
in at most $\delta_2 n$ steps, and when it terminates we are left
with a 4-graph $\cH'$ on $n'>(1-\delta_2)n$ vertices and  minimum
degree at least $d$.

Now suppose that we could prove that there is an edge of $\cH'$ that
lies in at least $(2-\e/2)(n'/4)^3$ copies of $P_2$. Since $\delta
\ll \e$, this is greater than $(2-\e)(n/4)^3$ and we are done. If on
the other hand $\cH'$ contains at least $2\delta n'^4$ copies of
$P_2$, then again this is at least $\delta n^4$ and we are done. So
if we could prove the result for $\cH'$ with $2\delta, \e/2$, then
we could prove the result for $\cH$ (with $\delta, \e$).
Consequently, we may assume that $\cH$ has minimum degree at least
$(3/32-o_{\delta}(1)){n \choose 3}=(1-\od1)(n/4)^3$.

If $\#P_2\ge \delta n^4$, then we are done so assume that
$\#P_2<\delta n^4=(\delta/n^3)n^7$.  Then by the Removal lemma, there is a set of at
most $\delta n^4$ edges of $\cH$ whose removal results in a
4-graph ${\cH}'$ with no copies of $P_2$. Since
$|{\cH}'|>t^4(n)-\delta n^4$, by Theorem \ref{P2stability},
we conclude that there is a partition of ${\cH}'$ (and also of
$\cH$) into four parts such that the number of edges intersecting some part in at least two points is $o_{\delta}(n^4)$. Now pick a partition $W \cup X \cup Y \cup Z$  of
$\cH$ that maximizes $e(W,X,Y,Z)$, the number of edges that intersect
each part. We know that $e(W,X,Y,Z)\ge t^4(n)-o_{\delta}(n^4)$, and an
easy calculation also shows that each of $W,X, Y, Z$ has size
$n/4\pm o_{\delta}(n)$.

Let $B$ be the set of edges of $\cH$ that intersect some part in at least two points. Let $G=\cH-B$ be the set of edges of $\cH$ that intersect each part.
  Let $M$ be the set of
4-tuples which intersect each part and are not edges of $\cH$. Then
$\cH-B \cup M= G \cup M$ is a 4-partite 4-graph with partition $W,X,Y, Z$,  so it has at most $t^4(n)$ edges. We
conclude that
\begin{equation} |M|<|B| <\odel(n^4),\end{equation}
in particular $B\ne \emptyset$.

\noindent{\bf Claim.}
For every vertex $a$ of $\cH$ we have $d_M(a)<\e_1(n/4)^3$, where $\e_1=\e/10^5$.

{\bf Proof of Claim.} Suppose for contradiction that $d_M(a)\ge \e_1(n/4)^3$ for some vertex $a$.  Then
$$(1-\od1)(n/4)^3\le d_{\cH}(a)=d_G(a)+d_B(a)\le (1+\od1)(n/4)^3-\e_1(n/4)^3+d_B(a).$$
We conclude that $d_B(a)\ge (\e_1-\od1)(n/4)^3 > (\e_1/2)(n/4)^3$. Let $L=L(a)$ be the set of triples $\{b,c,d\}$ such that $abcd \in B$.  So $|L|=d_B(a)>(\e_1/2)(n/4)^3$.
Partition $L=L_1 \cup L_2 \cup L_3$, where $L_i$ consists of those triples that intersect precisely $i$ parts.

Case 1: $|L_1| > (\e_1/10)(n/4)^3$.  Let us assume by symmetry that the number of triples $bcd \in L_1$ with $b,c,d \in W$ is at least $(\e_1/40)(n/4)^3$. For each choice of $(x,y,z) \in X \times Y \times Z$ with $a \ne x,y,z$, the three 4-tuples $bxyz, cxyz, abcd$ form a potential copy of $P_2$. The number of such choices of $(\{b,c,d\}, x,y,z)$ is at least $(1-\od1)(\e_1/40)(n/4)^6>\delta n^6$ so for at least half of these choices, one of the 4-tuples $bxyz, cxyz$ must be in $M$. Each of these 4-tuples is counted at most $n^2$ times, since $a$ is fixed. We obtain the contradiction
$(\e_1/100n^2)(n/4)^6<|M|=o_{\delta}(n^4)$.  This concludes the proof in
this case.

Case 2: $|L_2| > (3\e_1/10)(n/4)^3$.  Pick $bcd \in L_2$. There are $2\times{4 \choose 2}=12$ possibilities for the way the points $b,c,d$ are distributed within the parts.
Let us assume by symmetry that the number of triples $bcd \in L_2$ with $b,c\in W$, $d \in X$ is at least $(\e_1/40)(n/4)^3$.  Now proceed exactly as in the proof of Case 1.

Case 3: $|L_3| > (\e_1/10)(n/4)^3$. Assume wlog that $a \in W$. Pick $bcd \in L_3$. There are $3$ possibilities for the way the points $b,c,d$ are distributed within the parts (one point must be in $W$, the part containing $a$). Let us assume by symmetry that the number of triples $bcd \in L_3$ with $(b,c,d)\in W\times X \times Y$ is at least $(\e_1/30)(n/4)^3$. We may assume that
$d_{G}(a)\ge (\e_1/30)(n/4)^3$ for otherwise we can move $a$ to $Z$ and increase $e(W, X, Y, Z)$ thereby contradicting the choice of the partition. Now pick $bcd \in L_3$ as above and $(x,y,z) \in (X-\{c\}) \times (Y-\{d\}) \times Z$ with $axyz \in G$.  For each choice of
$(b,c,d,x,y,z)$ the three 4-tuples $abcd, axyz, bxyz$ form a copy of $P_2$.
The number of such choices of $(b,c,d,x,y,z)$ is at least $(\e_1/30)^2(n/4)^6>\delta n^6$ so for at least half of these choices, one of the 4-tuples $bxyz \in M$. Each of these 4-tuples is counted at most $n^2$ times, since $a$ is fixed. We obtain the contradiction
$(\e_1/30n)^2(n/4)^6<|M|=o_{\delta}(n^4)$.  This concludes the proof of
this case and the Claim.
\medskip

Partition $B=B_1 \cup B_2$, where $B_2$ consists of those edges of $B$ with exactly two points in one part, one point in a second part and one point in a third part (for example a $WWXY$ edge would be in $B_2$). Suppose that $B_1 \ne\emptyset$ and pick $e=abcd \in B_1$. Some two  points of $e$ must lie in the same part, so assume wlog that $a,b \in W$.

Let us first suppose that $c$ or $d$ is in $W$, say $c \in W$. For every $(x,y,z) \in X \times Y \times Z$
($x,y,z \ne d$), we get three potential copies of $P_2$ of the form
$wxyz, w'xyz, e$ where $w,w' \in\{a,b,c\}$. At least $(n/4)^3$ of these potential copies of $P_2$ contains a 4-tuple from $M$, otherwise we obtain $(2-\od1)(n/4)^3$
copies of $P_2$ containing $e$ and we are done.  Each such 4-tuple from $M$ is counted at most twice, so we obtain at least $(1/2)(n/4)^3$ 4-tuples from $M$ that contain some vertex in $\{a,b,c\}$. Consequently, there exists $w \in e$ with $d_M(w)\ge (1/6)(n/4)^3$ and this contradicts the Claim.

We may therefore assume that $a,b \in W$ and $c,d$ lie in the same part different from $W$, say $c,d \in X$. There are at least $(\e/3)(n/4)^3$ choices $(x,y,z) \in (X-\{c,d\}) \times Y \times Z$ with $vxyz \in M$ for some $v \in \{a,b\}$ or there are at least
$(\e/3)(n/4)^3$ choices $(w,y,z) \in (W-\{a,b\}) \times Y \times Z$ with $vxyz \in M$ for some $v \in \{c,d\}$.  This is because otherwise $e$ would lie in at least $(2-\od1-2\e/3)(n/4)^3>(2-\e)(n/4)^3$ copies of $P_2$.  In either case, we conclude that there exists $v \in e$ with $d_M(v)>(\e/6)(n/4)^3$ thus contradicting the Claim.

We conclude from the arguments above that $B_1=\emptyset$. Pick $e=abcd \in B_2$ and assume wlog that $a,b \in W, c \in X, d \in Y$.  For $(x,y,z) \in (X-\{c\}) \times (Y -\{d\})\times Z$, consider the following two potential copies of $P_2$:
$$e, bcdz, axyz \qquad \qquad e, acdz, bxyz.$$
The number of these potential copies of $P_2$ is twice the number of choices of $(x,y,z)$ and this is at least $(2-\od1)(n/4)^3$. At least $(\e/2)(n/4)^3$ of these potential $P_2$'s has a 4-tuple from $M$, otherwise we obtain at least $(2-\e)(n/4)^3$ copies of $P_2$ containing $e$ and we are done.  If for at least $(\e/4)(n/4)^3$ of these potential $P_2$'s, the 4-tuple from $M$ is of the form $wxyz, w \in \{a,b\}$
(i.e., the third in the lists), then there exists $v \in \{a,b\}$ with $d_M(v)\ge (\e/8)(n/4)^3>\e_1n^3$ thereby contradicting the Claim. So for at least $(\e/4)(n/4)^3$ of these potential $P_2$'s, the 4-tuple from $M$ is of the form $wcdz, w \in \{a,b\}$
(i.e., the second in the lists). Each such 4-tuple from $M$ is counted at most $|X||Y|$ times, so there are at least $(\e/4)(n/4)^3/(|X||Y|)>(\e/20)n$ 4-tuples from $M$ that intersect $e$ in three points.

Form the bipartite graph with parts $B=B_2$ and $M$ where
$e \in B$  is adjacent to  $f \in M$ if
$|e \cap f|=3$. We have shown above that each $e \in B$ has degree at least $(\e/20)n$. Since $|B|>|M|$, we conclude that there exists $f=abcd \in M$ (with $(a,b,c,d) \in W \times X \times Y \times Z$) that is adjacent to at least $(\e/20)n$ different $e$'s from $B$. Assume wlog that at least $(\e/80)n$ of these $e$'s contain $a,b,c$. We may also assume wlog that at least $(\e/240)n$ of these $e$'s have their fourth vertex in the same part as $a$, namely $W$.

Now for each $j=1, \ldots, (\e/240)n$, let $e_j=w_jabc \in B$ with $w_j \in W$. For every $(x,y,z) \in (X-\{b\}) \times (Y-\{c\}) \times Z$, consider the potential copy of $P_2$ given by $w_jxyz, axyz, e_j$.  The number of choices for $(j,x,y,z)$ is at least $(\e/240)n|X||Y||Z|>(4\e/10^5)n^4>2\delta n^4$.
If for at least half of these choices of $(j,x,y,z)$ the potential copy of $P_2$ is a real copy of $P_2$ in $\cH$, then  $\#P_2 \ge \delta n^4$ and we are done.
So we may assume that for at least half of the choices of $(j,x,y,z)$ (i.e. for at least $(2\e/10^5)n^4$ choices), the potential copy of $P_2$ referenced above has a 4-tuple $g \in M$. If at least half the time $a \not\in g$, then we obtain the contradiction
$(\e/10^5)n^4 \le |M| = o_{\delta}(n^4)$.  So  at least half the time $a \in g$.  Each such $g$ containing $a$ is counted at most $n$ times (once for each $w_j$), so we obtain $d_M(a)>(\e/10^5)n^3\ge \e_1 n^3$.  This contradicts the Claim and completes the proof of the theorem. \qed

\section{Counting $P_3$'s}

Recall that $c(n, P_3)=2(n/2)^3-\Theta(n^2)$.
Theorem \ref{books} for $l=3$ follows from the following result.

\begin{theorem} \label{P3}
For every $\e>0$ there exists $\delta>0$ and $n_0$ such that the
following holds for $n>n_0$. Every $n$ vertex 4-graph with
$d^4(n)+1$ edges contains either

$\bullet$ an edge that lies in at least $(2-\e)(n/2)^3$ copies of
$P_3$, or

$\bullet$ at least $\delta n^4
$ copies of $P_3$.
\end{theorem}

\bigskip
\noindent {\bf Proof of Theorem \ref{books} for $l=3$.} Remove $q-1$ edges
from $\cH$ and apply Theorem \ref{P3}.
If we find $\delta n^4$
copies of $P_3$, then since $q<\delta n$, the number of copies is
much larger than $q(1-\e)c(n, P_3)$ and we are done.
Consequently, we find
an edge $e_1$ in at least $(2-\e)(n/2)^3>(1-\e)c(n, P_3)$ copies of $P_3$. Now remove $q-2$ edges from $\cH-e_1$ and repeat this argument to obtain $e_2$.  In this way we obtain edges $e_1, \ldots, e_q$ as required.

The bound is sharp due to the following construction. Add a collection of $q$ pairwise disjoint edges within one part of $D^4(n)$.
It is easy to see that each added edge lies in
$2(n/2)^3+O(n^2)$ copies of $P_3$ and clearly no copy of $P_3$ contains two of
the new edges. Consequently, the copies of $P_3$ are counted exactly
once. \qed

We will need the following stability theorem for $P_3$ proved by F\"uredi-Pikhurko-Simonovits \cite{FPS}

\begin{theorem} {\bf ($P_3$ stability \cite{FPS})}\label{P3stability}
Let $\cH$ be a 4-graph with $n$ vertices and $d^4(n)-o(n^4)$ edges
that contains no copy of $P_3$. Then there is a partition of the
vertex set of $\cH$ into $X \cup Y$ so that the number of edges that
intersect some part in 0, 1, 3 or 4 points is $o(n^4)$. In other words, $\cH$ can
be obtained from $D^4(n)$ by adding and deleting a set of $o(n^4)$
edges.
\end{theorem}

 \bigskip

 \noindent{\bf Proof of Theorem \ref{P3}.}
 Given $\e$ let $0<\delta \ll \e$.
Write $o_{\delta}(1)$ for any function that approaches zero as
$\delta$ approaches zero and moreover, $o_{\delta}(1) \ll \e$. Let
$n$ be sufficiently large and let $\cH$ be an $n$ vertex 4-graph
with $d^4(n)+1$ edges.  Write $\#P_3$ for the number of copies of
$P_3$ in $\cH$.

As in the  proof of Theorem \ref{P2} (just replacing 3/32 by 3/8), we may assume that $\cH$ has minimum degree at
least  $d=(3/8)(1-\od1){n\choose 3}$

If $\#P_3\ge \delta n^4$, then we are done so assume that
$\#P_3<\delta n^4=(\delta/n^3)n^7$.  Then by the Removal lemma, there is a set of at
most $\delta n^4$ edges of $\cH$ whose removal results in a
4-graph ${\cH}'$ with no copies of $P_3$. Since
$|{\cH}'|>d^4(n)-\delta n^4$, by Theorem \ref{P3stability},
we conclude that there is a partition of ${\cH}'$ (and also of
$\cH$) into two parts such that the number of edges intersecting some part in 0,1,3, or 4 points is $o_{\delta}(n^4)$. Now pick a partition $X \cup Y$  of
$\cH$ that maximizes $e(X,Y)$, the number of edges that intersect
each part in two points. We know that $e(X,Y)\ge d^4(n)-o_{\delta}(n^4)$, and an
easy calculation also shows that each of $X, Y$ has size
$n/2\pm o_{\delta}(n)$.

Let $B$ be the set of edges of $\cH$ that intersect some part in 0,1,3 or 4 points. Let $G=\cH-B$ be the set of edges of $\cH$ that intersect each part in two points.
  Let $M$ be the set of
4-tuples which intersect each part in two points and are not edges of $\cH$. Then
$\cH-B \cup M= G \cup M$ is a  4-graph with $(2,2)$-partition $X \cup Y$,  so it has at most $d^4(n)$ edges. We
conclude that
\begin{equation} |M|<|B| <\odel(n^4),\end{equation}
in particular $B\ne \emptyset$.

\noindent{\bf Claim.}
For every vertex $a$ of $\cH$ we have $d_M(a)<\e_1(n/2)^3$, where $\e_1=\e/10^6$.

{\bf Proof of Claim.} Suppose for contradiction that $d_M(a)\ge \e_1(n/2)^3$ for some vertex $a$.  Then
$$(1/2-\od1)(n/2)^3\le d_{\cH}(a)=d_G(a)+d_B(a)\le (1/2+\od1)(n/2)^3-\e_1(n/2)^3+d_B(a).$$
We conclude that $d_B(a)\ge (\e_1-\od1)(n/2)^3 > (\e_1/2)(n/2)^3$. Let $L=L(a)$ be the set of triples $\{b,c,d\}$ such that $abcd \in B$.  So $|L|=d_B(a)>(\e_1/2)(n/2)^3$.  Assume wlog that $a \in X$.
Partition $L=L_{XXX} \cup L_{XXY} \cup L_{YYY}$, where $L_{X^iY^{3-i}}$ consists of those triples that intersect $X$ in precisely $i$ points (note that $L_{XYY}=\emptyset$ by definition of $B$).

Case 1: $|L_{XXX}| > (\e_1/6)(n/2)^3$ or $|L_{YYY}| > (\e_1/6)(n/2)^3$. Let us first assume that $|L_{XXX}| > (\e_1/6)(n/2)^3$. For each $bcd \in L_{XXX}$ with $e=abcd$ and  $(x,\{y,y'\}) \in (X-e) \times {Y\choose 2}$, the four 4-tuples $bxyy', cxyy', dxyy', e$ form a potential copy of $P_3$. The number of such choices of $(e, x,\{y,y'\})$ is at least $(\e_1/13)(n/2)^6>\delta n^6$ so for at least half of these choices, one of the 4-tuples $bxyy', cxyy', dxyy'$ must be in $M$. Each of these 4-tuples in $M$ is counted at most $|X|^2<n^2$ times, since $a$ is fixed. We obtain the contradiction
$(\e_1/13n^2)(n/2)^6<|M|=o_{\delta}(n^4)$.  If $|L_{YYY}| > (\e_1/6)(n/2)^3$, then the same proof works by replacing $(x, \{y,y'\})$ with $(\{x,x'\}, y) \in {X\choose 2} \times (Y-e)$. This concludes the proof in
this case.

Case 2: $|L_{XXY}| > (\e_1/6)(n/2)^3$.  We may assume that
$d_{G}(a)\ge |L_{XXY}|$ for otherwise we can move $a$ to $Y$ and increase $e(X, Y)$ thereby contradicting the choice of the partition. Pick $bcd \in L_{XXY}$ with $b,c \in X$ and $d \in Y$.  Consider $x,y,y'$ with $x \in X-e, y,y' \in Y-e$ and $axyy' \in G$.  For each choice of
$(e,x,\{y,y'\})$ the four 4-tuples $bcdy, bcdy', e, axyy'$ form a copy of $P_3$.
The number of such choices of $(e, x,\{y,y'\})$ is at least
$$d_G(a)|L_{XXY}| \ge |L_{XXY}|^2>2(\e_1/10^5)^2n^6>2\delta n^6$$ so for at least half of these choices, one of the 4-tuples $bcdy, bcdy' \in M$. Each of these 4-tuples of $M$ is counted at most $n^2$ times, since $a$ is fixed. We obtain the contradiction
$(\e_1/10^5)^2n^4<|M|=o_{\delta}(n^4)$.  This concludes the proof of
this case and the Claim.
\medskip

Partition $B=B_1 \cup B_2$, where $B_2$ consists of those edges of $B$ that intersect both parts in an odd number of points. Suppose that $B_1 \ne\emptyset$, pick $e=abcd \in B_1$ and assume wlog that $e \subset X$.
 For every $(x,\{y,y'\}) \in (X-e) \times {Y \choose 2}$, we get four potential copies of $P_3$ of the form
$w_1xyy', w_2xyy', w_3xyy', e$ where $\{w_1, w_2, w_3\} \subset e$ and $w_i\ne w_j$. At least $(\e/2)(n/2)^3$ of these potential copies of $P_3$ contains a 4-tuple from $M$, otherwise we obtain
\begin{equation} \label{star} 4(|X|-4){|Y| \choose 2}-(\e/2)(n/2)^3=(2-\od1-\e/2)(n/2)^3>(2-\e)(n/2)^3\end{equation}
copies of $P_3$ containing $e$ and we are done.  Each such 4-tuple from $M$ contains some point of $e$ so there exists $w \in e$
with  $d_M(w)\ge (\e/8)(n/2)^3>\e_1(n/2)^3$ and this contradicts the Claim.

We conclude the  $B_1=\emptyset$. Pick $e=abcd \in B_2$ and assume wlog that $a,b,c \in X, d \in Y$.  For $(x,\{y,y'\}) \in (X-e) \times {Y-e \choose 2}$, consider the following types of potential copies of $P_3$:

Type 1: $xyy'a, xyy'b, xyy'c, e$

Type 2: $e, abdy, abdy', xcyy'; \qquad e, acdy, acdy', xbyy'; \qquad e, bcdy, bcdy', xayy'$.

At least $(\e/2)(n/2)^3$ of these potential copies of $P_3$ contains a 4-tuple from $M$, otherwise we obtain at least $(2-\e)(n/2)^2$
copies of $P_3$ containing $e$ (as in (\ref{star})) and we are done. Suppose that at least half the time, the 4-tuple from $M$ is in one of the Type 1 copies, or the last 4-tuple in one of the type two copies (i.e., $xcyy', xbyy', xayy'$). Each such 4-tuple is counted at most twice, and so we obtain at least
$(\e/8)(n/2)^3$ 4-tuples of $M$ that intersect $e$. We conclude that there exists $w \in e$
with  $d_M(w)\ge (\e/32)(n/2)^3>\e_1(n/2)^3$ and this contradicts the Claim.

We may therefore assume that for at least $(\e/4)(n/2)^3$ of these potential copies of $P_3$, the 4-tuple from $M$ is one of the two middle ones of the Type 2 copies, and so it intersects $e$ in three points. Each such 4-tuple is counted at most $|X||Y|<(1+\od1)(n/2)^2$ times, so we obtain at least
$(\e/10)n$ 4-tuples from $M$ that intersect $e$ in three points. We have argued that for every $e \in B_2=B$, there are at least $(\e/10)n$ different $f \in M$ for which $|e \cap f|=3$. Since $|B|>|M|$, we conclude that there exists $f' \in M$ with at least $(\e/10)n$ different $e' \in B_2$ such that
$|e'\cap f'|=3$. At least $(\e/40)n$ of these $e'$'s intersect $f$ in the same three points. Consequently, we may assume wlog that there are $a,b \in X$, $d\in Y$ and $x_1, \ldots, x_{t} \in X$ with $t=(\e/40)n$ such that $e_i=abx_id \in B$.

Fix $i$, set $e=e_i$ and consider the Type 1 potential copies of $P_3$ referenced in the notation above with $c=x_i$, i.e., consider $xyy'a, xyy'b, xyy'c, e$.  Recall that there are at least $(1-\od1)|X|{|Y| \choose 2}>(1/20)n^3$ such copies. If at least $\e n^3$ of these potential copies of $P_3$ have a 4-tuple from $M$, then we find a vertex $w \in \{a,b,x_i\}$ with $d_M(w)>(\e/4)n^3>\e_1 n^3$ and this contradicts the Claim.
We conclude that  each $e_i=abcx_i$ lies in at least $(1/20-\e)n^3>\e n^3$ copies of $P_3$, and these copies are clearly distinct for distinct $i$. Altogether we therefore have $\#P_3 \ge t\e n^3=(\e^2/40)n^4>\delta n^4$ and we are done. \qed

\section{Counting $P_4$'s}
 Recall that
$c(n, P_4)=(4+o(1)){n/2 \choose 3}=\Theta(n^3)$.
Theorem \ref{books} for $l=4$ follows from the following result.

\begin{theorem} \label{P4}
For every $\e>0$ there exists $\delta>0$ and $n_0$ such that the
following holds for $n>n_0$. Every $n$ vertex 3-graph with
$b^4(n)+1$ edges contains either

$\bullet$ an edge that lies in at least $(1-\e)c(n, P_4)$ copies of
$P_4$, or

$\bullet$ at least $\delta n^4$ copies of $P_4$.
\end{theorem}

\bigskip
\noindent {\bf Proof of Theorem \ref{books} for $l=4$.}
Remove $q-1$ edges
from $\cH$ and apply Theorem \ref{P4}.
If we find $\delta n^4$
copies of $P_4$, then since $q<\delta n$, the number of copies is
much larger than $q(1-\e)c(n,P_4)$ and we are done.
Consequently, we find
an edge $e_1$ in at least $(1-\e)c(n,P_4)$ copies of $P_4$. Now remove $q-2$ edges from $\cH-e_1$ and repeat this argument to obtain $e_2$.  In this way we obtain edges $e_1, \ldots, e_q$ as required.

The result is asymptotically tight as we can add $q$ pairwise disjoint 4-tuples to $B^4(n)$, each intersecting both parts in two points.    \qed

\medskip

We need the following stability result proved in \cite{FMP}.

\begin{theorem} {\bf ($P_4$ stability \cite{FMP})}\label{P4stability}
Let $\cH$ be a 4-graph with $n$ vertices and $b^4(n)-o(n^4)$ edges
that contains no copy of $P_4$. Then there is a partition of the
vertex set of $\cH$ into $X \cup Y$ so that the number of edges that
intersect $X$ or $Y$ in an even number of points is $o(n^4)$. In other words, $\cH$ can
be obtained from $B^4(n)$ by adding and deleting a set of $o(n^4)$
edges.
\end{theorem}

 \bigskip

 \noindent{\bf Proof of Theorem \ref{P4}.}
 Given $\e$ let $0<\delta \ll \e$.
Write $o_{\delta}(1)$ for any function that approaches zero as
$\delta$ approaches zero and moreover, $o_{\delta}(1) \ll \e$. Let
$n$ be sufficiently large and let $\cH$ be an $n$ vertex 4-graph
with $b^4(n)+1$ edges.  Write $\#P_4$ for the number of copies of
$P_4$ in $\cH$.

As in the  proof of Theorem \ref{P2} (just replacing 3/32 by 1/2), we may assume that $\cH$ has minimum degree at
least  $d=(1/2-\od1){n\choose 3}$.

If $\#P_4\ge \delta n^4$, then we are done so assume that
$\#P_4<\delta n^4$.  Then by the Removal lemma, there is a set of at
most $\delta n^4$ edges of $\cH$ whose removal results in a
4-graph ${\cH}'$ with no copies of $P_4$. Since
$|{\cH}'|>b^4(n)-\delta n^4$, by Theorem \ref{P4stability},
we conclude that there is a partition of ${\cH}'$ (and also of
$\cH$) into two parts such that the number of edges intersecting some part in an even number of points is $o_{\delta}(n^4)$. Now pick a partition $X \cup Y$ of
$\cH$ that maximizes $e(X,Y)$, the number of edges that intersect
each part in an odd number of points. We know that $e(X,Y)\ge b^4(n)-o_{\delta}(n^4)$, and an
easy calculation also shows that each of $X, Y$ has size
$n/2\pm o_{\delta}(n)$.

Let $B$ be the set of edges of $\cH$ that intersect one (and therefore both) of $X, Y$ in an even number of points. Let $G=\cH-B$ be the set of edges of $\cH$ that intersect both $X, Y$ in an odd number of points.
  Let $M$ be the set of
4-tuples which intersect both parts in an odd number of points and are not edges of $\cH$. Then
$\cH-B \cup M= G \cup M$ is an odd 4-graph with partition $X,Y$,  so it has at most $b^4(n)$ edges. We
conclude that
\begin{equation} \label{mb} |M|<|B| <\odel(n^4),\end{equation}
in particular $B\ne \emptyset$. Let $B_{X^iY^{4-i}}$ ($\cH_{X^iY^{4-i}}$) be the set of edges in $B$ ($\cH$) with exactly $i$ points in $X$.  Let
$$\e_1=\min\{\e/200, \e^2/10^4, \e^3/10^3\}.$$

{\bf Claim.}
For every vertex $a$ of $\cH$ we have $d_M(a)<\e_1n^3$.

{\bf Proof of Claim.} Suppose for a contradiction that $d_M(a)>\e_1n^3$.  Since
$$(1/2-\od1){n \choose 3}\le d_{\cH}(a) = d_G(a)+d_B(a)\le \left((1/2+\od1){n \choose 3}-d_M(a)\right)+d_B(a)$$
we conclude that $d_B(a)>d_M(a)-\odel(n^3)>(\e_1/2)n^3$.  Assume wlog that $a \in X$. Then $d_B(a)=d_{B_{XXXX}}(a)+d_{B_{XXYY}}(a)$.

Case 1. $d_{B_{XXXX}}(a) \ge (\e_1/4)n^3$

We may assume that $d_G(a)\ge d_{B_{XXXX}}(a)$ otherwise moving $a$ from $X$ to $Y$
increases $e(X,Y)$ and contradicts the choice of the partition.
Pick $e=abcd \in B_{XXXX}$ and $f=auvw \in G$. The number of choices for $(e,f)$ is at least $d_{B_{XXXX}}(a)^2\ge (\e_1/4)^4n^6$.  For each such $(e,f)$, consider the five 4-tuples
$$uvwb, \quad uvwc, \quad uvwd, \quad f, \quad e.$$
This gives at least $(\e_1^2/16)n^6>2\delta n^6$ potential copies of $P_4$ so at least $(\e_1^2/32)n^6$ of them have a 4-tuple not in $\cH$.  Since $e,f \in \cH$, the absent 4-tuple is of the form $uvwz$ where $z \in \{b,c,d\}$. Notice that $uvwz \in M$ since $auvw \in G$ and $a,z \in X$. Each such 4-tuple of $M$ is counted at most $3n^2$ times, as there are at most three choices for $|e \cap f|$ and ${|X| \choose 2}<n^2$ choices for the two vertices of $e-f-\{a\}$. This yields the contradiction
$(\e_1^2/32)n^6/(3n^2) <|M| =\odel(n^4)$.

Case 2. $d_{B_{XXYY}}(a) \ge (\e_1/4)n^3$

Suppose that $d_{\cH_{XYYY}}(a)\ge (\e_1/2)n^3$.  Then pick $e=abcd \in B_{XXYY}$ and $f=auvw \in \cH_{XYYY}$ (so $u,v,w \in Y$) with $|e \cap f|=1$.  The number of such pairs $(e,f)$ is at least $(\e_1^2/10)n^6$.  For each such $(e,f)$, consider the potential copy of $P_4$ given by
$$bcdu, \quad bcdv, \quad bcdw, \quad e, \quad f.$$
Since $(\e_1^2/10)n^6>2\delta n^6$ at least half of them have a 4-tuple from $M$. Each such 4-tuple is counted at most $3n^2$ times, so we obtain the contradiction
$(\e_1^2/20)n^6/(3n^2)\le |M| =\odel(n^4)$. We may therefore assume that

(i) $d_{\cH_{XYYY}}(a)< (\e_1/2)n^3$ and

(ii) $d_{B_{XXXX}}(a) < (\e_1/4)n^3$.

Define $L(a)=\{bcd: abcd \in \cH\}$, so $|L(a)|=d_{\cH}(a)\ge (1/2-\od1){n \choose 3}$. Consider the partition
$$L(a)=L_{XXX} \cup L_{XXY} \cup L_{XYY} \cup L_{YYY},$$
where the subscripts have the obvious meaning. Then (i) and (ii) translate to
$$|L_{YYY}|+|L_{XXX}| <\e_1n^3.$$
For $(u,v) \in X \times Y$, let $d_L(uv)$ be the number of $w$ such that $uvw \in L(a)$. Then
\begin{align}
\sum_{(u,v) \in X \times Y} d_L(u,v)&=2(|L_{XXY}|+|L_{XYY}|) \notag \\
&=2(|L(a)|-|L_{YYY}|-|L_{XXX}|)\notag \\
&>2(|L(a)|-\e_1 n^3)\notag \\
&\ge 2(1/2-\od1
){n \choose 3}-2\e_1n^3 \notag \\
&\ge (1-13\e_1 -\od1){n \choose 3}.\notag
\end{align}
Consequently, there exists $(b,c) \in X \times Y$ such that
$$d_L(bc) >\frac{(1-13\e_1 -\od1){n \choose 3}}{(1/4+\od1)n^2}>(2/3-9\e_1)n.$$
We conclude that there exists $S \subset X, T \subset Y$ such that
$$\min\{|S|, |T|\}\ge (2/3 -9\e_1-1/2-\od1)n>(1/6-10\e_1)n$$
  and $abcd \in \cH$ for every $d \in S \cup T$.  Now pick $s_1, s_2, s_3 \in S$ and $t \in T$ and consider the potential $P_4$
\begin{equation} \label{r}abcs_1, \quad abcs_2, \quad abcs_3, \quad abct, \quad s_1s_2s_3t.\end{equation}
The number of choice for $(\{s_1, s_2, s_3\}, t)$ is at least ${|S| \choose 3}|T|>10^{-4}n^4$. If for at least half of these choices of $(\{s_1, s_2, s_3\}, t)$ we get a copy of $P_4$ in $\cH$ as shown above, then $\#P_4>(1/2)10^{-4}n^4>\delta n^4$, a contradiction. So for at least half of the choices, one of the 4-tuples in (\ref{r}) is not in $\cH$. By definition of $S$ and $T$, the first four are in $\cH$, so the last one is in $M$. This is counted exactly once, so we obtain the contradiction $(1/2)10^{-4}n^4<|M|=\odel(n^4)$.
This completes the proof of the Claim.

\smallskip
Partition $B=B_1 \cup B_{XXYY}$ where
$$B_1=B_{XXXX} \cup B_{YYYY}$$

{\bf Case 1.} $|B_1| \ge \e|B|$.

Pick $e=abcd \in B_1$, and assume wlog that $e \in B_{XXXX}$. Let
$e' \subset e$ with  $|e'|=3$.  Assume wlog that $e'=bcd$. Let $\{y_1, y_2, y_3\} \in {Y \choose 3}$ and consider the five 4-tuples
$$bcdy_1, \quad bcdy_2, \quad bcdy_3, \quad e,\quad ay_1y_2y_3.$$
These 4-tuples from a potential copy of $P_4$.  The number of choices for $(e', \{y_1, y_2, y_3\})$ is at least $4(1-\od1){n/2 \choose 3}$. For at least $2\e{n/2 \choose 3}$ of these choices, one of the 4-tuples above must not be in $\cH$, otherwise
$\#P_4\ge 4(1-\od1){n/2 \choose 3}-2\e{n/2 \choose 3}>(1-\e)c(n, P_4)$ and we are done. If for at least $\e{n/2 \choose 3}$ of these choices, the missing 4-tuple is the last one in the list, then we obtain $d_M(x)>(\e/4){n/2 \choose 3}>\e_1n^3$ for some $x \in e$ (since $\e_1 \le \e/200$). This contradicts the Claim. We may therefore assume that for at least $\e{n/2 \choose 3}$ of these choices, the 4-tuple from $M$ has exactly three points in $e$.
Each such 4-tuple is counted at most ${|Y| \choose 2}$ times giving at least $(\e/7)n$ 4-tuples from $M$ with three points in $e$.

We have argued above that for every $e \in B_1$, there are more than $(\e/7)n$ different $f \in M$ with $|e \cap f|=3$. Form the bipartite graph with parts $B_1, M$, where $e \in B_1$ is adjacent to $f \in M$ if $|e \cap f|=3$.  Then each vertex in $B_1$ has degree more than $(\e/7)n$, and since $|B_1|\ge \e|B| > \e|M|$, we conclude that there exists $f \in M$ that is adjacent to more than $|B_1|(\e/7)n/|M|> (\e^2/7)n$ different $e \in B_1$. Consequently, there exist $a,b,c$ such that $d_{B_1}(abc)> (\e^2/7)n$.  Assume wlog that $a,b,c \in X$.

For each choice of $d$ with $e=abcd \in B_1$ and $\{y_1, y_2, y_3\} \in {Y \choose 3}$ five 4-tuples $y_1y_2y_3x$ where $x \in e$ together with $e$ form a potential copy of $P_4$. The number of choices for $(d, \{y_1, y_2, y_3\})$ is at least $d_{B_1}(abc){|Y| \choose 3}>6\e_1n^4>2\delta n^4$ (since $\e_1\le \e^2/10^4$).  If for at least half of them, we get a copy of $P_4$ in $\cH$, then $\#P_4 >\delta n^4$ and we are done.  So for at least $3\e_1n^4$ of the choices, one of the five 4-tuples is not in $\cH$.  If for at least $\e_1 n^4$ choices the missing 4-tuple is of the form  $y_1y_2y_3d$, then we obtain the contradiction
$\e_1 n^4\le |M|=\odel(n^4)$. So for at least $2\e_1 n^4$ choices the missing 4-tuple is of the form  $y_1y_2y_3x$, $x\ne d$. Each such missing 4-tuple is counted at most $|X|<n$ times. We conclude that there exists $x \in e$ with $d_M(x)>2\e_1 n^4/n>\e_1n^3$.  This contradicts the Claim and completes the proof in this case.
\medskip

{\bf Case 2.} $|B_1| < \e|B|$.

In this case we have $|B_{XXYY}| \ge (1-\e)|B|$.  Partition $B_{XXYY}=B_2 \cup B_3$
where
$$B_2=\{e \in B_{XXYY}: d_M(e')>(1-\e)(n/2) \hbox{ for every $e' \subset e$ with $|e'|=3$}\}.$$
Suppose that $|B_2|\ge (1-\e)|B_{XXYY}|$. Then we count 4-tuples of $M$ from sets in $B_2$.  For each set in $B_2$, there are four choices for $e' \subset e$ with $|e'|=3$, and given $e'$, there are  $(1-\e)(n/2)$ 4-tuples of $M$ containing $e'$. Each 4-tuple from $M$ is counted at most $3\max\{|X|, |Y|\}$ times. This gives the contradiction
$$|M| \ge \frac{4(1-\e)(n/2)|B_2|}{3\max\{|X|, |Y|\}}>\frac{4(1-2\e)|B_2|}{3}
\ge \frac{4(1-2\e)(1-\e)^2|B|}{3}>|B|>|M|.$$
We may therefore suppose that $|B_3|>\e |B_{XXYY}|>(\e/2)|B|$.  We may also assume that no edge of $B$ lies in at least $(1-\e)c(n, P_4)$ copies of $P_4$, otherwise we are done. Using this observation we conclude that  we have at least $(\e/4)n$ 4-tuples in $M$. To see this we pick an edge $e \in B$ and consider potential copies of $P_4$
containing $e$. We know that at least $(\e/2)c(n, P_4)$ of these potential copies
have a 4-tuple from $M$, for otherwise $e$ lies in at least $(1-\e)c(n, P_4)$ copies of $P_4$. Each such 4-tuple is counted at most $\max\{{|X| \choose 2}, {|Y| \choose 2}\}$ times.
So we may assume that $$|B|>|M|>(\e/4)n.$$
 Now pick an edge  $e =abcd\in B_3$.  By definition of $B_3$, there exists  $e'=bcd \subset e$ with $d_M(bcd)\le (1-\e) n/2$. Assume wlog that $b \in X, c,d \in Y$.
Then there is a set $Y' \subset Y$ such that $bcdy \in \cH$ for every $y \in Y'$ and
$$|Y'| \ge |Y|-2-d_M(bcd)\ge (1-\od1)(n/2)-(1-\e)(n/2)>(\e/4)n.$$ By the Claim and $\e_1\le \e^3/10^3$, we know that the number of $\{y_1, y_2, y_3\} \in {Y'\choose 3}$ with $ay_1y_2y_3 \in \cH$ is at least
$${|Y'| \choose 3}-d_M(a)>{(\e/4)n \choose 3}-\e_1n^3>2\e_1 n^3-\e_1n^3=\e_1n^3.$$
Each such $\{y_1, y_2, y_3\} \in {Y'\choose 3}$ yields the $P_4$ given by
$$bcdy_1,\quad bcdy_2,\quad bcdy_3,\quad e,\quad ay_1y_2y_3.$$
We have argued above that for each  $e \in B_3$ there are at least $\e_1n^3$ copies of $P_4$ containing $e$. Each such copy of $P_4$ contains a unique edge of $B_3$.  Consequently, we obtain
$$\#P_4 \ge |B_3|(\e_1n^3)>(\e/2)|B|(\e_1n^3)>(\e^2\e_1/8)n^4>\delta n^4.$$
This contradiction completes the proof of the theorem. \qed

\section{Counting Expanded triangles}

Recall that $c(n, C_3)=3(n/2)^2+\Theta(n)$.
Theorem \ref{C3q} follows from the following result.

\begin{theorem} \label{C3}
For every $\e>0$ there exists $\delta>0$ and $n_0$ such that the
following holds for $n>n_0$. Every $n$ vertex 4-graph with
$b^4(n)+1$ edges contains either

$\bullet$ an edge that lies in at least $(3-\e)(n/2)^2$ copies of
$C_3$, or

$\bullet$ at least $\delta n^4$ copies of $C_3$.
\end{theorem}

\bigskip
\noindent {\bf Proof of Theorem \ref{C3q}.}
Remove $q-1$ edges
from $\cH$ and apply Theorem \ref{C3}.
If we find $\delta n^4$
copies of $C_3$, then since $q<\delta n^2$, the number of copies is
much larger than $q(1-\e)c(n, C_3)$ and we are done.
Consequently, we find
an edge $e_1$ in at least $(3-\e)(n/2)^2>(1-\e)c(n, C_3)$ copies of $C_3$. Now remove $q-2$ edges from $\cH-e_1$ and repeat this argument to obtain $e_2$.  In this way we obtain edges $e_1, \ldots, e_q$ as required.

 Sharpness
follows by the following construction: Add a collection of $q$
4-tuples to $B^4(n)$ within one of the parts (say $X$) such that  every two 4-tuples have at most
one point in common. It is well-known that such quadruple-systems
exist of size $\delta n^2$ (in fact such Steiner
 systems also
exist for an appropriate congruence class of $n$).  It is easy to see that each added 4-tuple lies in at most
$3(n/2)^2$ copies of $C_3$, since there are three ways to partition the edge into two disjoint pairs, and for each of these ways, there are at most $(n/2)^2$ copies of $C_3$ using this partition. Moreover, no two added edges lie in a copy of $C_3$ since they share at most one point.
Consequently, the number of copies of $C_3$ is at most
$3q(n/2)^2$.   \qed

We need the following stability result proved by Keevash and Sudakov \cite{KSFr}.

\begin{theorem} {\bf ($C_3$ stability \cite{KSFr})}\label{C3stability}
Let $\cH$ be a 4-graph with $n$ vertices and $b^4(n)-o(n^4)$ edges
that contains no copy of $C_3$. Then there is a partition of the
vertex set of $\cH$ into $X \cup Y$ so that the number of edges that
intersect $X$ or $Y$ in an even number of points is $o(n^4)$. In other words, $\cH$ can
be obtained from $B^4(n)$ by adding and deleting a set of $o(n^4)$
edges.
\end{theorem}

 \bigskip

 \noindent{\bf Proof of Theorem \ref{C3}.}
 Given $\e$ let $0<\delta \ll \e$.
Write $o_{\delta}(1)$ for any function that approaches zero as
$\delta$ approaches zero and moreover, $o_{\delta}(1) \ll \e$. Let
$n$ be sufficiently large and let $\cH$ be an $n$ vertex 4-graph
with $b^4(n)+1$ edges.  Write $\#C_3$ for the number of copies of
$C_3$ in $\cH$.

As in the  proof of Theorem \ref{P2} (just replacing 3/32 by 1/2), we may assume that $\cH$ has minimum degree at
least  $d=(1/2-\od1){n\choose 3}$.

If $\#C_3\ge \delta n^4$, then we are done so assume that
$\#C_3<\delta n^4=(\delta/n^2)n^6$.  Then by the Removal lemma, there is a set of at
most $\delta n^4$ edges of $\cH$ whose removal results in a
4-graph ${\cH}'$ with no copies of $C_3$. Since
$|{\cH}'|>b^4(n)-\delta n^4$, by Theorem \ref{C3stability},
we conclude that there is a partition of ${\cH}'$ (and also of
$\cH$) into two parts such that the number of edges intersecting some part in an even number of points is $o_{\delta}(n^4)$. Now pick a partition $X \cup Y$ of
$\cH$ that maximizes $e(X,Y)$, the number of edges that intersect
each part in an odd number of points. We know that $e(X,Y)\ge b^4(n)-o_{\delta}(n^4)$, and an
easy calculation also shows that each of $X, Y$ has size
$n/2\pm o_{\delta}(n)$.

Let $B$ be the set of edges of $\cH$ that intersect one (and therefore both) of $X, Y$ in an even number of points. Let $G=\cH-B$ be the set of edges of $\cH$ that intersect both $X, Y$ in an odd number of points.
  Let $M$ be the set of
4-tuples which intersect both parts in an odd number of points and are not edges of $\cH$. Then
$\cH-B \cup M= G \cup M$ is an odd 4-graph with partition $X,Y$,  so it has at most $b^4(n)$ edges. We
conclude that
\begin{equation} \label{m<b} |M|<|B| <\odel(n^4),\end{equation}
in particular $B\ne \emptyset$.

Given vertices $a,b$ and hypergraph $F$, write $d_F(ab)$ for the number of edges of $F$ containing both $a$ and $b$.  The rest of the proof has many similarities (modulo technical changes) to the proof of the exact Tur\'an result for $C_3$ in \cite{KSFr}.
Let $\e_1=\e/10^5$.
\medskip

{\bf Claim 1.} For every two vertices $a,b \in V:=X \cup Y$, either $d_G(ab)<\e_1 n^2$ or $d_B(ab)<\e_1 n^2$.

Proof.  Suppose, for contradiction, that both  $d_G(ab)$ and $d_B(ab)$ are at least $\e_1n^2$.
Pick $e \in B$ and $f \in G$ with $e \cap f=\{a,b\}$.  Note that in all cases $g_{e,f}=e \cup f -\{a,b\} \in M \cup G$, i.e., $g_{e,f}$ intersects both parts in an odd number of points. The number of such pairs $e,f$ is at least $(\e_1n^2)^2/2>2\delta n^4$  (the  factor of 2 is to ensure that $e \cap f=\{a,b\}$). If at least $\delta n^4$ of these pairs form a copy of $C_3$, then we are done, so we may assume that at least $(\e_1n^2)^2/4 $ of these pairs satisfy $g_{e,f} \in M$. This contradicts (\ref{m<b}) and completes the proof of the Claim.

{\bf Claim 2.} $d_B(v)<(\e/10^3) n^3$ for every $v \in V$.

Proof. Let us fix a vertex $v \in V$, $\e'=\e/10^3>24\e_1$, and assume for contradiction that $d_B(v)\ge \e' n^3$. Call vertex $w \in V-\{v\}$ good if $d_B(vw)<\e_1 n^2$, otherwise say that $w$ is bad. Claim 1 implies that if $w$ is bad, then $d_G(vw)<\e_1n^2$. Moreover, the number of bad vertices is at least $\e' n$ for otherwise we obtain the contradiction
$$d_B(v) \le \sum_{w\, bad} d_B(vw) +\sum_{w\, good} d_B(vw) < \e' n {n \choose 2} + n \e_1n^2<\e' n^3.$$
Next we observe that $d_G(v) \ge d_B(v)$ for otherwise we could move $v$ to the other part and contradict the choice of $X,Y$. This implies that $d_G(v) \ge (1/4-o_{\delta}(1)){n \choose 3}$. If the number of good vertices is less than $n/18$, then  $$d_G(v) \le \sum_{w\, good} d_G(vw) +\sum_{w\, bad} d_G(vw) <{n\over 18}{n \choose 2}+n(\e_1 n^2)<\left({1\over6}+7\e_1\right){n \choose 3}.$$
This contradicts the lower bound on $d_G(v)$. We may therefore assume that the number of good vertices is $\alpha n$, where
\begin{equation} \label{ep} 1/18\le \alpha \le 1-\e'.\end{equation}
Write $d_{\cH}(v)=d_G(v)+d_B(v)$ and let us estimate these two terms separately.
The number of edges of $G$ containing $v$ and a bad vertex is at most $((1-\alpha)n+1)\e_1 n^2<\e_1 n^3$. The  number of edges of $G$ containing $v$ and no bad vertex is at most
$${\alpha n(\alpha n-1)\over 6}\left({1\over 2}+o_{\delta}(1)\right)n\le \left({\alpha^2\over 2} +o_{\delta}(1)\right){n \choose 3}.$$
The bound above is obtained by picking two good vertices which then restricts the edge being counted to one of the parts. This procedure counts each edge six times.
We conclude that $d_G(v) <(\alpha^2/2+6\e_1){n \choose 3}$.

The number of edges of $B$ containing $v$ and a good vertex is at most $\e_1 n^3$.  Using a similar argument to that used above, the number of edges of $B$ containing $v$ and no good vertex is at most $((1-\alpha)^2/2+6\e_1){n \choose 3}$. We conclude that
$$d_{\cH}(v)\le \left({\alpha^2 +(1-\alpha)^2\over 2} +12\e_1\right){n \choose 3}.$$
Using (\ref{ep}) and $\e_1<\e'/24$, we observe that
 $$\frac{\alpha^2+(1-\alpha)^2}{2}+12\e_1=\frac12 +\alpha^2-\alpha+12\e_1<\frac12 +(1-\e')^2-(1-\e')+12\e_1<\frac12 -\frac{\e'}{2}.$$
 Consequently,
 $d_{\cH}(v)<(1/2-\e'/2){n \choose 3}$. This contradicts the fact that $d_{\cH}(v)\ge (1/2-o_{\delta}(1)){n \choose 3}$ and completes the proof of the Claim.

\medskip

The rest of the proof is devoted to showing that $d_B(v)\ge (\e/10^3) n^3$ for some vertex $v$ and this contradicts Claim 2. Note that for every edge $e \in B$, there are at least $(3-o_{\delta}(1))(n/2)^2$ copies of $C_3$ containing $e$ where the other two edges in the copy are in $G$.  Indeed, this is why $c(n, C_3)=(3+o(1))(n/2)^2$.
This requires some case analysis, for example, if $e=\{a,b,c,d\} \subset X$, then for every choice of $(x,y) \in (X-e) \times Y$, and for every partition of $e$ into two disjoint pairs $p,q$, the three edges $e, p \cup \{x,y\}, q \cup \{x,y\}$ form a copy of $C_3$ and $p \cup \{x,y\}, q \cup \{x,y\} \in G$. The number of such copies is therefore the number of $(x,y)$ times the number of pairs $p,q$ and this is $(3-o_{\delta}(1))(n/2)^2$.
The case $a,b \in X, c,d \in Y$ is similar except that the argument further breaks into two cases depending on the choice of $p,q$. We omit these details.

{\bf Claim 3.} There is a pair of vertices $a,b$ with
$d_B(ab)>(\e/48) n^2$

Proof. For every $e \in B$, at least $(\e/2)(n/2)^2$ of the copies of $C_3$ using $e, f, g$ with $f,g \in G\cup M$ must have at $f\in M$ or $g \in M$.  Otherwise, there are at least $(3-o_{\delta}(1)-\e/2)(n/2)^2>(3-\e)(n/2)^2$
 copies of $C_3$ containing $e$ and we are done. The edge in $M$ is counted precisely once, since a copy of $C_3$ is uniquely determined by two of its edges. We conclude that for each $e \in B$, there are at least $(\e/2)(n/2)^2$ edges of $M$ that intersect $e$ is exactly two points. Now form a bipartite graph with parts $B,M$ where $e \in B$ is adjacent to $f \in M$ if $|e \cap f|=2$.  Since $|M|<|B|$, and each $e \in B$ is adjacent to at least  $(\e/2)(n/2)^2$ different $f \in M$, we conclude that there exists $f \in M$ that is adjacent to more than $ (\e/2)(n/2)^2$ different $e \in B$. At least 1/6 of these $e$ intersect $f$ in the same pair of points $a,b$.
 Consequently, $d_B(ab)>  (\e/12)(n/2)^2 =(\e/48)n^2$ and the Claim is proved.
 \medskip

 Let us fix $a,b$ from Claim 3. For each edge $e=abcd \in B$, there are at least $n^2/5$ pairs $r,s \in V$ such that the three sets $e, acrs, bdrs$ form a copy of $C_3$. By Claim 3, the number of such potential copies of $C_3$ is at least $(\e/240)n^4>2\delta n^4$, so for at least half of them, either $acrs \in M$ or $bdrs \in M$.  Each such 4-tuple of $M$ is counted at most $n$ times, since $a,b, r,s$ are fixed. This gives us at least $(\e/480)n^3$ 4-tuples of $M$ containing either $a$ or $b$.  At least $(\e/960)n^3$ must contain the same point, say $a$. Consequently, $d_M(a)\ge (\e/960)n^3$. We know that $d_{\cH}(a)\ge (1/2-o_{\delta}(1)){n \choose 3}$, and the above argument shows that $d_G(a) \le (1/2+o_{\delta}(1)-\e/960){n \choose 3}$.  We conclude that $d_B(a)>(\e/10^3)n^3$ which contradicts Claim 2 and completes the proof.  \qed


\end{document}